\documentclass[preprint,11pt]{article}
 \usepackage{amsthm,amsmath,amssymb}

\newtheorem{thm}{Theorem}

\newtheorem{Rem}{Remark}
\usepackage{float}
\usepackage{graphicx}
\usepackage{dsfont}
\usepackage[table]{xcolor}

\usepackage{enumerate}

\newtheorem{corollary}{Corollary}

\title{A Consistent Estimator for  Skewness of Partial Sums of Dependent Data}

\author{Masoud M Nasari and Mohamedou  Ould-Haye}
\date{\small  School of Mathematics and Statistics.
\\
Carleton University, 1125 Colonel By Dr. Ottawa, ON, Canada, K1S 5B6}

\begin{document}
\maketitle
\begin{abstract}
We introduce an estimation method for the scaled {\it skewness coefficient of the sample mean}  of short and long memory linear processes. This method can  be extended to estimate higher moments  such as  curtosis coefficient of the sample mean. Also a general result on computing all asymptotic moments of partial sums is obtained, allowing in particular a much easier derivation of some existing central limit theorems for linear processes. The introduced skewness estimator  provides a tool to empirically  examine  the error  of the central limit theorem for long and short memory linear processes. We also show that, for both  short and long memory linear processes,  the   skewness coefficient of the sample mean converges to zero at the same rate  as in the i.i.d. case.
\end{abstract}

\section{Introduction}\label{Introduction}
Skewness is a characteristic of a distribution which  is often used to measure its departure  from symmetry. The skewness of the marginal distribution of a set of data is a different quantity from that of an estimator which is built as an aggregation of the data. In particular, the skewness of $\overline{X}$, which is the sample mean of $n$, $n \geq 1$ i.i.d. observations with a finite third moment, is of the same order as  $1/\sqrt{n}$ times the skewness of the marginal distribution. In other words, the marginal skewness of i.i.d. is the same as $\sqrt{n}$ times the skewness of  $\overline{X}$.  As a result, in the i.i.d. case, having an empirical estimator for the marginal skewness means having an estimator for the skewness of the sample mean. In the case of stationary and dependent  observations, the relation between the skewness of the sample mean and that of the marginal distribution of the data is no longer as  straightforward as that of  the i.i.d. data. For stationary dependent data, unlike in the i.i.d. case, estimating the skewness of the marginal distribution and estimating that of the sample mean do not amount to the same task. Although, the marginal skewness of the stationary data has been investigated and contributions have been made to the area (cf. e.g. Bai and Ng (2005), Grigoletto and Lisi (2009) and  references therein) less so has been done for their sample mean. In this paper, a consistent empirical estimator for the third cumulant of partial sum of $n$ stationary linear processes is introduced.  This estimator is then combined with  some existing consistent estimators for the variance of partial sums of linear processes to  obtain an estimator for the skewness of partial sums of linear processes. The results in this paper also show that the rate at which the skewness of partial sums of long and short memory linear processes vanishes is, counter-intuitively for long memory, $\sqrt{n}$ and that long memory effect tends to make the sampling distribution more and more quickly symmetrical. The rest of the paper is organized as follows: in section 2  we establish and discuss (the impact of) a general result (Theorem 1) on asymptotically computing all higher order moments of partial sums of linear processes. This result is interesting on its own as, given its short proof,  it opens the door to  reproving some existing central limit theorems without using the  machinery of stochastic integrals or cumulant-based methods to prove asymptotic normality. In section 2 we  also give a second important result which is the skewness estimation (Theorem 2) where we essentially develop a method to estimate the rightly normalized third moment of partial sum. This method can be  extended to estimate all higher moments. In section 4 we discuss some Monte Carlo studies to illustrate the skewness estimation. Section 5 is devoted to proofs.

\section{Skewness of partial sums of linear processes}\label{Skewness of partial sums of linear processes}
For throughout use in this paper we let   $X_1,\cdots,X_n$,   be a stationary sample of size $n$, $n\geq 1$, with the linear representation
%Let us consider a stationary process $(X_t)$ with a linear representation

\begin{equation}\label{Definition linear processes}
X_t=\mu+\sum_{i=0}^\infty a_i \varepsilon_{t-i}, ~~ t\geq 1,
\end{equation}
where $\mu \in \mathds{R}$, $a_i$ are squarely summable and $\varepsilon_i$, $i\geq 0$, are i.i.d. white noise with variance $\sigma^2$, where $0<\sigma^2<\infty$. Furthermore,   we assume that $\varepsilon_1$ has a finite third moment.
The validity of some of the results in this paper will require that $\varepsilon_1$ has a finite sixth moment.
\\
The linear process  (\ref{Definition linear processes}),  when  $\sum_{i=0}^{\infty} |a_{i}|< \infty$ is said to be  of short memory. For long memory linear processes, we adopt the conventional definition in which $X_t$ is called long memory if, as $i \to \infty$ and for some constant $c(d)>0$ (depending on $d$),  $a_i \sim c(d) i^{d-1}$,  with $0< d <1/2$, where, here and also  throughout this paper, $\ell_n \sim  \lambda_n$ stands for the asymptotic  equivalence    $\ell_n=\lambda_n(1+o(1))$.
\\
The parameter $d$ is  the memory parameter whose large values (i.e., closer to 1/2) indicate a stronger dependence between the observations on the linearly structured temporal  series $X_t$ as in (\ref{Definition linear processes}).
\\
To unify our notation we denote a short memory  linear process $X_t$, as in (\ref{Definition linear processes}), with    $d=0$.
\\
We note that linear processes of the form (\ref{Definition linear processes}) englobe the well known fractionally autoregressive  integrated moving average FARIMA($p,d,q)$ processes of the form
$$
\phi(B)(1-B)^d(X_t-\mu)=\theta(B)\epsilon_t
$$ where $\phi(B)=1-\phi_1B-\cdots-\phi_pB^p$ and $\theta(b)=1+\theta_1B+\cdots+\theta_qB^q$ are respectively the $p$ order and $q$ order   autoregressive and  moving average polynomials, and $B$ is the backshift operator, namely $BX_t=X_{t-1}$. From corollary 3.1 of Kokoszka and Taqqu (1995), such processes have an explicit  linear representation (\ref{Definition linear processes}) with
$$
a_i\sim\frac{\theta(1)}{\phi(1)}\frac{i^{d-1}}{\Gamma(d)},\qquad0<d<1/2,\qquad \textrm{as }i\to\infty,
$$
and when $d=0$, i.e. just ARMA($p,q$), the coefficients $a_i$ decrease exponentially fast towards zero. In all what follows and without restriction of generality, we assume that $\mu=0$. Let $\mathbb{E}(\epsilon_1^3)=\eta$. Also, let $S_n=X_1+\cdots X_n$.
Before presenting the skewness estimation, we first give a general result on asymptotically computing all higher order moments of partial sums of linear processes according to wether the  higher order is even or odd.
\begin{thm}\label{Theroem rate of S_3}
Assume that $X_t$ satisfies (\ref{Definition linear processes}) with $0\le d<1/2$, and that $\mathbb{E}[\epsilon_1^k]$ is finite. Let
\begin{equation}\label{md}
m(d)=\begin{cases}\sum_{i=0}^\infty a_i,&\textrm{if }d=0,\\\\
\frac{c(d)}{d},&\textrm{if }0<d<1/2.
\end{cases}
\end{equation}
Then as $n\to\infty$,
$$
n^{-p(1+2d)}\mathbb{E}\left[S_n^k\right]\to
(m(d))^k\sigma^k\left(\frac{k!}{2^p(p!)}\right)\left(\frac{1}{1+2d}+\int_0^\infty\left((1+x)^d-x^d\right)^2dx\right)^p
$$
if $k=2p$, and
\begin{align*}
&\sqrt{n}\,\,n^{-\frac{k}{2}(1+2d)}\mathbb{E}\left[S_n^k\right]\to(m(d))^k\eta\sigma^{k-3}\binom{3+2\ell}{3}\left(\frac{(2\ell)!}{2^\ell(\ell!)}\right)\times\\
&
\left(\frac{1}{1+3d}+\int_0^\infty\left((1+x)^d-x^d\right)^3dx\right)
\left(\frac{1}{1+2d}+\int_0^\infty\left((1+x)^d-x^d\right)^2dx\right)^\ell
\end{align*}
if $k=3+2\ell$, $\ell\ge0$.
\end{thm}
\begin{Rem}
The previous theorem provides another (much easier method of moments) proof of  some Central limit theorems for linear processes (see for example Davydov 1970) when all moments exist. Actually, from one side, it is well know that if $Z=\mathcal{N}(0,1)$
then for all integer $p\ge0$,
$$\mathbb{E}(Z^{2p})=\frac{(2p)!}{2^p(p!)},\qquad\textrm{and }\mathbb{E}(Z^{2p+1})=0,
$$
and from another side, it is straightforward from the theorem above to see that
$$
\mathbb{E}\left[\left(\frac{S_n}{\sqrt{\textrm{Var}(S_n)}}\right)^k\right]
\to\begin{cases}\frac{(2p)!}{2^p(p!)},&\textrm{if }k=2p,\\\\0,&\textrm{if }k=2p+1.\end{cases}
$$
We also note that the previous theorem goes deeper than showing that odd moments of normalized partial sums converge to zero and gives an explicit rate of $\sqrt{n}$ of this convergence regardless of the type of the memory whether it be short or long. \end{Rem}
\begin{Rem} We note that for all FARIMA processes, the coefficient $m(d)$ in (\ref{md}) is continuous in $d$ as $d\to0$. Actually, and as mentioned prior to  the previous theorem, for FARIMA($p,d,q$), (and noting that
$\Gamma(d)\sim 1/d$ as $d\to0$), $c(d)=\theta(1)/(\phi(1)\Gamma(d))\sim(\theta(1))/(\phi(1))d$ so that $m(d)\to m(0)$ as $d\to0$,   since FARIMA($p,0,q)$ is just ARMA$(p,q)$, for which the coefficients $a_i$ in (\ref{Definition linear processes}) satisfy
$$
\sum_{k=0}^\infty a_k=\frac{\theta(1)}{\phi(1)},
$$
since
$$
\sum_{k=0}^\infty a_kz^k=\frac{\theta(z)}{\phi(z)},\qquad\textrm{for all }\vert z\vert\le1.
$$
\end{Rem}
We now define  a measure of   skewness  for the partial sum of $X_1,\cdots,X_n$, i.e., the first $n\geq 1$ observations of the linear process $X_t$, as follows:

\begin{equation}\label{skewness}
  \beta_n:= \frac{\mathbb{E}(S_n^3)}{\left[\textrm{Var}(S_n)\right]^{3/2} }.
\end{equation}
%For the ease of the notation and without loss of generality   from this stage henceforth we  assume  that    $\mu=E(X_t)=0$.
%\\
Note that $\beta_n=0$ for all $n$ for   linear processes  $X_t$, as   in (\ref{Definition linear processes}), with symmetric innovations, as we will have in this case $\mathbb{E}(\epsilon_1^3)=0$, but in general, due to the CLT, we will always have $\beta_n\to0$ as $n\to\infty$.
\\
The rate at which the skewness $\beta_n$ of the partial sums of a linear process vanishes is given in the following corollary which is an immediate consequence of Theorem \ref{Theroem rate of S_3} by evaluating the second  and the third moments of $S_n$.
\begin{corollary}\label{Corollary rate of beta}
Let $X_1,\cdots,X_n$ be the first $n$ terms of the  linear process (\ref{Definition linear processes}) with $0 \leq d <1/2$ fixed.  If $\eta$ is finite then, as $n\to\infty$,
\begin{equation}\label{skew limit}
\sqrt{n}\beta_n\to
k(d):=\frac{\eta}{\sigma^3}\frac{\frac{1}{1+3d}+\int_0^\infty\left((1+x)^d-x^d\right)^3dx}{
\left(\frac{1}{1+2d}+\int_0^\infty\left((1+x)^d-x^d\right)^2dx\right)^{3/2}}.
\end{equation}
\end{corollary}
\begin{Rem}
From the previous corollary  one can see that the skewness of sums of  short and long memory linear processes asymptotically vanishes at the same convergence rate $\sqrt{n}$ as that of   sums of i.i.d. data. The effect of the short or long range dependence appears only in terms of $d$ in the limiting constant $k(d)$ in (\ref{skew limit}) which, surprisingly enough, (in absolute value)  is  decreasing  in $d\in[0,1/2)$, with  $k(0)=\eta/\sigma^3$ and $k(.5)=0$.  While short memory effect (such as in ARMA models) tends to asymptotically disappear, long memory effect tends to make normalized partial sums more and more quickly symmetrical. This is mainly due to the fact that, at fixed $n$  and  as $d\to.5$, the variance of the  partial sum becomes infinite, while  its third moment remains bounded. Of course these findings have to be considered with some caution as we still do not know how fast (as function of $n$ and $d$) convergence (\ref{skew limit}) is taking place.
\end{Rem}
\section{A consistent  estimator for skewness of partial sums of linear processes}\label{A consistent  estimator for skewness of partial sums of linear processes}
To construct   a consistent  estimator for $\beta_n$, one approach is to construct  separate consistent estimators for $\mathbb{E}\left(S_n\right)^3$  and $\textrm{Var}(S_n)$. In this section, we estimate only the former  as    a  consistent estimator   for the variance of the partial sums of short and long memory linear process is obtained in Abadir \textit{et al} (2009). More precisely, Abadir {\it et al}  used the well known fact that
$$
n^{-1-2d}\textrm{Var}\left(\sum_{k=1}^nX_k\right)=n^{-2d}\left(\gamma(0)+2\sum_{h=1}^n\left(1-\frac{h}{n}\right)\gamma(h)\right)\to v(d)
$$
where $v(d)>0$ is the limiting variance and  $\gamma(h)=\textrm{Cov}(X_1,X_{1+h})$, to construct the long run variance estimator from the sample covariance function
$$
\hat\gamma(h)=\frac{1}{n}\sum_{j=1}^{n-h}(X_j-\bar X)(X_{j+h}-\bar X),
$$
in the form of
$$
q_0^{-2d }  \left(\bar{\gamma}_0  +2 \sum_{h=1}^{q_0}  \left(1-\frac{h}{q_0}\right) \bar{\gamma}_h\right)\overset{P}{\to} v(d),
$$
where $q_0\to\infty$ and $q_0=o(n)$.
\\
To   estimate  the limiting normalized third moment of the partial sum  (for $0\le d<1/2$)
\begin{equation}\label{Definition of S_3}
\mathcal{S}_3(d):=n^{-1-3d}\mathbb{E}(S_n^3)\to\frac{\eta}{\sigma^3}m(d)\left(\frac{1}{1+3d}+\int_0^\infty\left((1+x)^d-x^d\right)^3dx\right),
\end{equation}
we define the estimator
\begin{equation}\label{Definition of bar S_3}
\overline{\mathcal{S}}_3(d):=q_1^{-3d}\overline{\Delta}(0)+3q_2^{-3d}\sum_{h=1}^{q_2}\left(1-\frac{h}{q_2}\right)\overline{\Delta}(h)+
6q_3^{-3d}\sum_{h=1}^{q_3-1}\sum_{h'=1}^{q_3-h}\left(1-\frac{h+h'}{q_3}\right)\overline{\Delta}(h,h'),
\end{equation} where, as $n\to\infty$, $q_i\to\infty$, for $i=1,2,3$, with  $q_2=o(n)$,  $q_3=o(n^{1/2})$,  and $\overline{\Delta}(h)$ and $\overline{\Delta}(h,h')$   are third order sample covariances  defined as

\begin{eqnarray*}
\overline{\Delta}(h)&:=&\frac{1}{n}\sum_{j=1}^{n-h}\left((X_j-\bar X)^2(X_{j+h}-\bar X)+(X_j-\bar X)(X_{j+h}-\bar X)^2\right),\\
\overline{\Delta}(h,h')&:=& \frac{1}{n}\sum_{j=1}^{n-h-h'}(X_j-\bar X)(X_{j+h}-\bar X)(X_{j+h+h'}-\bar X).
\end{eqnarray*}
The idea of having such estimator stems from the fact that we can break $\mathcal{S}_3(d)$ into similar third order theoretical covariances (according to the cases where all indices are equal, or two are equal, or all are different)
\begin{equation}\label{third}
\mathcal{S}_3(d)=n^{-3d}\left(\Delta(0)+3\sum_{h=1}^n\left(1-\frac{h}{n}\right)\Delta(h)+6\sum_{h=1}^{n-1}\sum_{h'=1}^{n-h}\left(1-\frac{h+h'}{n}\right)
\Delta(h,h')\right),
\end{equation}
where
$$
\Delta(h)=\mathbb{E}(X_1X^2_{1+h})+\mathbb{E}(X_1^2X_{1+h}),\qquad \Delta(h,h')=\mathbb{E}(X_1X_{1+h}X_{1+h+h'}).
$$
In practice we would use $\mathcal{S}_3(\hat d)$ where $\hat d$ is a consistent estimator for $d$.
The bandwidth numerical sequence  $q_n=q$ is  such that, as $n \to \infty$,    $q \to \infty$ and  $q^2=o(n)$.
\\
The following theorem is  the second  main result of this paper which shows that $\overline{\mathcal{S}}_3(d)$ is   a consistent estimator for $\mathcal{S}_3(d)$, the normalized third moment of  partial sum of linear processes.

\begin{thm}\label{theorem2}
If $\mathbb{E}(\varepsilon_1^6)<\infty$ then, as $n, q_i \to\infty$ in such a way that $q_i^2=o(n)$, for $i=1,2,3$ used in the definition (\ref{Definition of bar S_3}), we have \begin{equation*}
\overline{\mathcal{S}}_3(d)-\mathcal{S}_3(d)=o_P(1).
\end{equation*}
\end{thm}
A direct consequence of Theorem \ref{theorem2} and Abadir {\it et al} yields a consistent  estimator for the scaled skewness $\sqrt{n}\beta_n$, where $\beta_n$ is as in (\ref{skewness}). This consistency is stated in the following corollary.
\begin{corollary}
Assume that  $\mathbb{E}(\varepsilon_1^6)<\infty$. As  $n, q_0,q_i, \to\infty$ in such a way that $q_i^2=o(n)$, for $i=1,2,3$, and $q_0=o(n)$, we have
\begin{equation*}
\hat k(d)\to k(d),
\end{equation*}
where $k(d)$ is given in (\ref{skew limit}) and
\begin{equation}\label{skewness estimator}
\hat k(d):=\frac{\overline{\mathcal{S}}_3}{\left( q_0^{-2d }  (\bar{\gamma}_0  +2 \sum_{h=1}^{q_0}  (1-h/q_0) \bar{\gamma}_h)   \right)^{3/2}}.
\end{equation}\end{corollary}

\section{Numerical illustrations}
The results in Tables 1 and 2 illustrate the
performance of the empirical estimator $\hat k(d)$ above,  in
estimating  $k(d)$. In Tables 1 and 2, for each sample size, 2000
replications of an  ARMA(1,1) process $X_t=\phi X_{t-1}+ \theta \varepsilon_{t-1}+\varepsilon_t$ were produced. In Tables 3 and 4 the same number of replications were produced from a fractionally
integrated process $X_t=(1-B)^{-d} \varepsilon_t$, where $d$ is the memory
parameter, $B$ is the back shift  operator. In both cases, the innovations $\varepsilon_t$  are i.i.d. exponentially distributed with mean 1.
In all the
tables, the following  empirical mean square error  is computed:
\begin{equation*}
  \widehat{MSE}\big(\hat k(d)\big)=\frac{\sum_{b=1}^{2000}
\Big( \hat k_b(d)  -   k(d)\Big)^2 }{2000},
\end{equation*}
where $\hat k_b(d) $ is  computed based on the $b$th replication
of the linear process, as indicated in each  table.
\begin{table*}[!htb]
    %\caption{Global caption}
    \begin{minipage}{.5\linewidth}
      \caption{\protect$X_t=0.5X_{t-1}+0.5\varepsilon_{t-1}+\varepsilon_t$}
      \centering
        \begin{tabular}{c c c c c c c c   }
\hline
& & & & & &  \\
& & $n$ &  &  $\widehat{MSE}(\hat k(d))$ &   &  &\\
 \hline
& & 200  & &1.075& & & \\
& & 1000 & &0.575& & &   \\
& & 5000 & &0.298& & &    \\
\hline
\end{tabular}
    \end{minipage}%
    \begin{minipage}{.5\linewidth}
      \centering
        \caption{\protect $X_t=-0.5X_{t-1}-0.5\varepsilon_{t-1}+\varepsilon_t$}
        \begin{tabular}{ c c c c c c c c }
\hline
 & & & & & & & \\
& & $n$ & &  $\widehat{MSE}(\hat k(d))$ & & &  \\
 \hline
& & 200 &  &1.800& & & \\
& & 1000&  &0.586& & &   \\
& & 5000&  &0.171& & &      \\
\hline
\end{tabular}
    \end{minipage}
\end{table*}

\begin{table*}[!htb]
    %\caption{Global caption}
    \begin{minipage}{.5\linewidth}
      \caption{ARIMA(0,0.2,0)}
      \centering
        \begin{tabular}{c c c c c c c c   }
\hline
& & & & & & & \\
& & $n$ &  & $\widehat{MSE}(\hat k(d))$ &  &  &  \\
 \hline
& & 200 &  &0.374& & & \\
& & 1000&  &0.113& & &    \\
& & 5000&  &0.048& & &     \\
\hline
\end{tabular}
    \end{minipage}%
        \begin{minipage}{.5\linewidth}
      \centering
        \caption{\protect ARIMA(0,0.4,0)}
        \begin{tabular}{c c c c c  c  c c  }
\hline
& & & & & & &  \\
& & $n$ &  &$\widehat{MSE}(\hat k(d))$ & & &    \\
 \hline
& & 200&  &0.121& & &\\
& &1000&  &0.027& & &    \\
& &5000&  &0.023& & &     \\
\hline
\end{tabular}
    \end{minipage}
\end{table*}
In the preceding tables, the limiting values of the normalized skewness coefficient for the two ARMA(1,1) and ARIMA(0,0.2,0) and ARIMA(0,0.4,0), when innovations $\epsilon_t$ are exponentially distributed (i.e. with skewness coefficient 2), are respectively 2 and 1.7 and .675, as we immediately  see from (\ref{skew limit}) that $k(0)=2$ and we can  numerically compute $k(0.2)=1.7$ and $k(.4)=.675$.  It is worth reiterating  that for all short memory linear processes, the limiting normalized skewness coefficient $k(0)$ is equal to the skewness of the innovations $\epsilon_t$.
\begin{Rem}
It was observed empirically that for long memory processes with all memory parameters $0<d<1/2$, the choice $q_1=q_2=\lceil n^{0.2} \rceil, ~ q_3=\lceil n^{0.1} \rceil, ~ q_4=\lceil n^{0.5-d} \rceil$, in (\ref{Definition of bar S_3}) and (\ref{skewness estimator}) yields relatively good estimates for $k(d)$. For short memory processes, the choice $q_2=q_3= q_4=\lceil n^{0.33} \rceil$ is deemed to be a good choice. The preceding choices were implemented in  the numerical illustrations of Tables 1-4.
\end{Rem}
\section{Proofs}
\subsection{ Proof of Theorem \ref{Theroem rate of S_3}}
Assume that $0<d<1/2$. With the assumption that $\mu=0$ and writing (\ref{Definition linear processes}) as
$$
X_t=\sum_{j=-\infty}^ta_{t-j}\epsilon_j
$$
we get
\begin{eqnarray}\label{odd-even}\lefteqn{
\mathbb{E}\left[\left(\sum_{t=1}^nX_t\right)^k\right]
=\mathbb{E}\left[\left(\sum_{j=-\infty}^n\left(\sum_{t=\max(j,1)}^na_{t-j}\right)\epsilon_j\right)^k\right]}\nonumber\\
&&=\sum_{j_1=\infty}^n\cdots\sum_{j_k=\infty}^n\left(\sum_{t=\max(j_1,1)}^na_{t-j_1}\right)\cdots\left(\sum_{t=\max(j_k,1)}^na_{t-j_k}\right)\mathbb{E}\left(\epsilon_{j_1}\cdots\epsilon_{j_k}
\right).
\end{eqnarray}
We show that in (\ref{odd-even}), in the case $k=2p$, only the indices $j_1,\ldots,j_k$ that are equal two by two will be the leading terms, and in case $k=3+2p$, only the cases corresponding two three of them are equal and the rest are equal to by two will lead. Consider first the case $k=2p$. The two by two equal indices give
$$
\left(\frac{k!}{2^p(p!)}\right)\sigma^k\left[\left(\sum_{j=-\infty}^n\left(\sum_{t=\max(j,1)}^na_{t-j}\right)^2\right)^p-B_n\right],
$$
where
$$
B_n=\sum_{s=1}^{p-1}\,\,\sum_{\alpha_1+\cdots+\alpha_s=p}\,\,\sum_{j_1=-\infty}^n\cdots\sum_{j_s=-\infty}^n\left(\sum_{t=\max(j_1,1)}^na_{t-j_1}\right)^{2\alpha_1}\cdots
\left(\sum_{t=\max(j_s,1)}^na_{t-j_s}\right)^{2\alpha_s}
$$
with $\alpha_1,\ldots,\alpha_s\ge1$. The coefficient $k!/(2^p(p!))$ corresponds to the number of two-by-two configurations which is the number  of pairings of $k=2p$ individuals.
Now we have (distinguishing between $j\ge1$ and $j\le0$), as $n\to\infty$,
\begin{eqnarray*}\lefteqn{
\sigma^k\left(\sum_{j=-\infty}^n\left(\sum_{t=\max(j,1)}^na_{t-j}\right)^2\right)^p}\\
&&=\sigma^k\left(\sum_{j=1}^n\left(\sum_{t=j}^na_{t-j}\right)^2+\sum_{j=0}^\infty\left(\sum_{t=1}^na_{t+j}\right)^2\right)^p\\
&&=\sigma^k\left(\sum_{j=1}^n\left(\sum_{i=0}^ja_{i}\right)^2+\sum_{j=1}^\infty\left(\sum_{t=1}^na_{t+j}\right)^2\right)^p\\
&&\sim\sigma^k\left(\left(\frac{c(d)}{d}\right)^2\frac{1}{2d+1}n^{2d+1}+\sum_{j=1}^\infty\left(\frac{c(d)}{d}(n+j)^d-(j+1)^d\right)^2\right)^p\\
&&=\sigma^k\left(\frac{c(d)}{d}\right)^k\left(\frac{1}{2d+1}n^{2d+1}+\sum_{j=1}^\infty\left((n+j)^d-(j+1)^d\right)^2\right)^p\\
&&\sim\sigma^k\left(\frac{c(d)}{d}\right)^kn^{(2d+1)p}\left(\frac{1}{2d+1}+\int_0^\infty\left((x+1)^d-x^d\right)^2dx\right)^p
\end{eqnarray*}
For the case of short memory ($d=0$), where the coefficients $a_i$ are summable, we have
$$
\sum_{j=1}^\infty\left(\sum_{t=1}^na_{t+j}\right)^2=\sum_{t=1}^n\sum_{s=1}^n\sum_{j=1}^\infty a_{t+j}a_{s+j}=
\sum_{t=1}^n\sum_{j=1}^\infty a_{t+j}\sum_{s=1}^na_{s+j}
=\sum_{t=1}^no(1)=o(n),
$$
and
$$
\sum_{j=1}^n\left(\sum_{i=0}^ja_{i}\right)^2\sim n\left(\sum_{i=0}^\infty a_i\right)^2.
$$
Hence in both cases of short and long memory, and when $k=2p$, we get that
\begin{eqnarray*}\lefteqn{
n^{-p-2dp}\left(\sum_{j=-\infty}^n\left(\sum_{t=\max(j,1)}^na_{t-j}\right)^2\right)^p}\\
&&\to (m(d))^k\left(\frac{1}{2d+1}+\int_0^\infty\left((x+1)^d-x^d\right)^2dx\right)^p.
\end{eqnarray*}
Similar calculations show that
$$
B_n\sim\sum_{s=1}^{p-1}(m(d))^kn^{2sd+s}\,\,\sum_{\alpha_1+\cdots+\alpha_s=p}\,\,\prod_{j=1}^s\left(\frac{1}{2d\alpha_j+1}+\int_0^\infty\left((x+1)^d-x^d\right)^{2\alpha_j}
\right)=o\left(n^{2pd+p}\right).
$$
Now for the indices that are not two by two equal, and similarly to $B_n$, they will yield quantities of the form
$$
\mathbb{E}(\epsilon_1^{\beta_1})\cdots\mathbb{E}(\epsilon_1^{\beta_s})n^{kd+s}(m(d))^k
\prod_{j=1}^s\left(\frac{1}{d\beta_j+1}+\int_0^\infty\left((x+1)^d-x^d\right)^{\beta_j}
\right)
$$
where $s=1,\ldots,p-1$, and $\beta_j\ge2$ and $\beta_1+\cdots+\beta_s=k.$ Clearly each quantity above is $o(n^{2pd+p})$. This completes the first part of the theorem (i.e. when $k$ is even). Now we consider $k=3+2\ell$, with $\ell\ge0$, we proceed in a similar way as above for the $2\ell$ part as we must have three indices (among $j_1,\ldots,j_k$) equal and the remaining $2p$ indices must be two by two equal in order to build  the leading term. I.e. when $k=3+2\ell$, we have $\alpha(k)=:\binom{k}{3}(2p)!/(2^p(p!))$ ways of selecting three indices (to be equal) and  paring the remaining $2p$, so that (\ref{odd-even}) is asymptotically equivalent to
\begin{align*}
&\alpha(k)\sum_{j_0=-\infty}^n\sum_{j_1=-\infty}^n\cdots\sum_{j_\ell=-\infty}^n\left(\sum_{t=\max(j_0,1)}^na_{t-j_0}\right)^3\left(\sum_{t=\max(j_1,1)}^na_{t-j_1}\right)^2
\cdots\left(\sum_{t=\max(j_\ell,1)}^na_{t-j_\ell}\right)^2
\mathbb{E}(\epsilon_1^3)\left(\mathbb{E}(\epsilon_1^2)\right)^\ell \\
&\sim\alpha(k)\eta\sigma^{k-3}(m(d))^kn^{3d+1}n^{2d\ell+\ell}\left(\frac{1}{3d+1}+\int_0^\infty\left((x+1)^d-x^d\right)^3dx\right)
\left(\frac{1}{2d+1}+\int_0^\infty\left((x+1)^d-x^d\right)^2dx\right)^\ell &\\
&=\alpha(k)\eta\sigma^{k-3}(m(d))^kn^{kd+\frac{k-1}{2}},
\end{align*}
which completes the proof of the second part of the theorem.
\subsection{Proof of Theorem \ref{theorem2}}
We first    introduce the following notations.
$$
\hat\Delta(h)=\frac{1}{n}\sum_{j=1}^{n-h}\left[(X_j-\mu)^2(X_{j+h}-\mu)+(X_j-\mu)(X_{j+h}-\mu)^2\right],
$$
and
$$
\hat\Delta(h,h')=\frac{1}{n}\sum_{j=1}^{n-h}(X_j-\mu)(X_{j+h}-\mu)(X_{j+h+h'}-\mu).
$$
Let
$$
\widehat{\mathcal{S}}_3=q_1^{-3d}\hat\Delta(0)+3q_2^{-3d}\sum_{h=1}^{q_2}\left(1-\frac{h}{q_2}\right)\hat\Delta(h)+
6q_3^{-3d}\sum_{h=1}^{q_3-1}\sum_{h'=1}^{q_3-h}\left(1-\frac{h+h'}{q_3}\right)\hat\Delta(h,h').
$$

If $\mathbb{E}(\varepsilon_1^6)<\infty$ then, as $n\to\infty$, we have
\begin{equation}\label{lem1}
\overline{\mathcal{S}}_3-\widehat{\mathcal{S}}_3=o_P(1).
\end{equation}
Equation (\ref{lem1}) follows from the law of large numbers for the sample mean of linear processes.

In view of  (\ref{lem1}), the conclusion of  this theorem is equivalent to

\begin{equation}\label{lem5}
\widehat{\mathcal{S}}_3-\mathcal{S}_3=o_P(1).
\end{equation}
To prove  (\ref{lem5}), we assume that $\mu=0$ and,  to make notations lighter, we will drop the coefficients $(1-h/n), (1-h/q_i),(1-(h+h')/n),(1-(h+h')/q_i)$  and replace all of them by their limit (as $n,q_i\to\infty$)  1. Also,   from the proof of Theorem 1, we have   for $0<d<1/2$,
$$
\Delta(h)=\eta\left(\sum_{i=0}^\infty a_i^2a_{i+h}+\sum_{i=0}^\infty a_ia^2_{i+h}\right) \sim \mathbb{E}(X_1^2X_{1+h}).
$$
By virtue of the preceding asymptotic equivalence,  we take
$$
\Delta(h)=\mathbb{E}(X_1^2X_{1+h}),
$$
and accordingly  we will also consider that
\begin{equation}\label{asymp1}
\hat\Delta(h)=\frac{1}{n}\sum_{j=1}^{n-h}X_j^2X_{j+h}.
\end{equation}
Clearly $\hat\Delta(0)\to\Delta(0)$, so, we just need to deal with  the two remaining terms in $
\widehat{\mathcal{S}}_3$. That is we need to show that
\begin{equation}\label{eq1}
q_2^{-d}\left(\sum_{h=1}^{q_2}\hat\Delta(h)-\sum_{h=1}^{q_2}\Delta(h)\right)=o_P(1)
\end{equation}
and that
\begin{equation}\label{eq2}
q_3^{-3d}\left(\sum_{h=1}^{q_3}\sum_{h'=1}^{q_3-h}\hat\Delta(h,h')-\sum_{h=1}^{q_3}\sum_{h'=1}^{q_3-h}\Delta(h,h')\right)=o_P(1).
\end{equation}
To prove (\ref{eq1}), without loss of generality and just to avoid the use of absolute value at different places, we will assume that $\Delta(h)\ge0$ and $a_j\ge0$ and $a_j$ are nonincreasing. Then, in view of (\ref{asymp1}), we have
$$
q_2^{-d}\left\vert\sum_{h=1}^{q_2}\mathbb{E}\left(\hat\Delta(h)-\Delta(h)\right)\right\vert=
q_2^{-d}\sum_{h=1}^{q_2}\frac{h}{n}\Delta(h)\le\frac{q_2}{n}q_2^{-d}\sum_{h=1}^{q_2}\Delta(h)=O(q_2/n).
$$
Consequently,
\begin{eqnarray}\label{right}
q_2^{-2d}\mathbb{E}\left(\sum_{h=1}^{q_2}\left(\hat\Delta(h)-\Delta(h)\right)^2\right)&\sim& q_2^{-2d}\textrm{Var}
\left(\sum_{h=1}^{q_2}\hat\Delta(h)\right)\nonumber\\
&\sim& q_2^{-2d}\mathbb{E}\left(\left(\sum_{h=1}^{q_2}\hat\Delta(h)\right)^2\right)-q_2^{-2d}\left(\sum_{h=1}^{q_2}\Delta(h)\right)^2.
\end{eqnarray}
We first deal with the second term of the right hand side of  (\ref{right}) as follows.
\begin{equation*}%\label{gammah}
{q_2}^{-d}\sum_{h=1}^{q_2}\Delta(h)={q_2}^{-d}\sum_{h=1}^{q_2}\mathbb{E}\left[\left(\sum_{i=0}^\infty a_i\varepsilon_{1-i}\right)^2\sum_{j=0}^\infty a_j\varepsilon_{1+h-j}\right]=\eta {q_2}^{-d}\sum_{h=1}^{q_2}\sum_{k=0}^\infty a_k^2a_{k+h}\to\eta c_1(d).
\end{equation*}
For the first term of the right hand side of  (\ref{right}) we write
\begin{align}\label{sumgamma}
&&{q_2}^{-2d}\mathbb{E}\left(\left(\sum_{h=1}^{q_2}\hat\Delta(h)\right)^2\right)
={q_2}^{-2d}\sum_{h=1}^{q_2}\sum_{h'=1}^{q_2}\mathbb{E}(\hat\Delta(h)\hat\Delta(h'))
=\frac{{q_2}^{-2d}}{n^2}\sum_{h=1}^{q_2}\sum_{h'=1}^{q_2}\sum_{i=1}^{n-h}\sum_{j=1}^{n-h'}\mathbb{E}\left(X_i^2X_{i+h}X_j^2X_{j+h}\right)\nonumber\\
%&&=\frac{{q_2}^{-2d}}{n^2}\sum_{h=1}^{q_2}\sum_{h'=1}^{q_2}\sum_{i=1}^{n-h}\sum_{j=1}^{n-h'}\sum_{s=0}^\infty
&&=\frac{{q_2}^{-2d}}{n^2}\sum_{h=1}^{q_2}\sum_{h'=1}^{q_2}\sum_{i=1}^{n-h}\sum_{j=1}^{n-h'}\sum_{s=0}^\infty
%%% I ADDED Eta^2 here %%%%%
\sum_{t=0}^\infty\sum_{f=0}^\infty\sum_{s'=0}^\infty\sum_{t'=0}^\infty\sum_{f'=0}^\infty
a_sa_ta_fa_{s'}a_{t'}a_{f'}\mathbb{E}\left[\varepsilon_{i-s}\varepsilon_{i-t}\varepsilon_{i+h-f}\varepsilon_{j-s'}\varepsilon_{j-t'}\varepsilon_{j+h'-f'}\right].\nonumber\\
&&
\end{align}
We now show that only the configuration corresponding to the three by three equality of indices  $i-s=i-t=i+h-f\neq j-s'=j-t'=j+h'-f'$
will contribute to the limit as $n$ (and therefore ${q_2}$) tends to infinity, in the equation (\ref{sumgamma}). This configuration corresponds to the following quantity,with $\eta=\mathbb{E}(\varepsilon_1^3$), the subtracted term corresponds to taking off the diagonal elements  with  $t'=t+k$ where $k=j-i$:
\begin{eqnarray*}
\eta^2{q_2}^{-2d}\sum_{h=1}^{q_2}\sum_{h'=1}^{q_2}\left(\sum_{t=0}^\infty a_t^2a_{t+h}\right)\left(\sum_{t=0}^\infty a_t^2a_{t+h'}\right)
-\frac{\eta^2{q_2}^{-2d}}{n}\sum_{h=1}^{q_2}\sum_{h'=1}^{q_2}\sum_{t=0}^\infty\sum_{k=-n}^n a_t^2a_{t+h}a^2_{t+k}a_{t+k+h'}.
\end{eqnarray*}
%The last term being nonnegative and clearly bounded by
The preceding term is nonnegative and bounded above by
$$
\frac{{q_2}^{-2d}}{n}\sum_{h=1}^{q_2}\sum_{h'=1}^{q_2}\Delta(h)\Delta(h')=O(1/n).
$$
%Hence we see that the considered  three by three configuration above results in the limit we are seeking
As a result, with the configuration $i-s=i-t=i+h-f\neq j-s'=j-t'=j+h'-f'$, the term on the right hand side of the last equation in (\ref{sumgamma}),  as $n,{q_2} \rightarrow \infty$ such that ${q_2}=o(n^{1/2})$,  converges to the following nonzero limit:
%%%% I ADDED THIS PART HERE %%%%%
$$
\underset{q\to\infty}{\lim} {q_2}^{-2d}\sum_{h=1}^{q_2}\sum_{h'=1}^{q_2}\Delta(h)\Delta(h')=\underset{{q_2}\to\infty}{\lim}\left({q_2}^{-d}\sum_{h=1}^{q_2}\Delta(h)\right)^2.
$$
All other configurations of the indices  of the summations on the right hand side of (\ref{sumgamma}),
result in asymptotically negligible terms, as $n,{q_2} \rightarrow \infty$.  In what follows we show the negligibility of selected configurations of the indices, noting that the remaining cases of each configuration can be treated in the same way.\\
%In what follows we pick an example of each of the remaining configurations in the sum (\ref{sumgamma}) to show that they are asymptotically negligible, noting that the remaining cases of each configuration can be treated in the same way.
%Among the remaining cases of the three by three configuration
We now consider the  configurations  $i-s=i-t=j-s'\neq
i+h-f=j-t'=j+h'-f'$, which is equivalent to
$s=t$, $s'=t+k$, with $k=j-i$, $f=t'+h+k$, $f'=t'+h'$,  and  $i-s=i-t=i+h-f= j-s'=j-t'=j+h'-f'$.  Defining %$M=(\mathbb{E}(\varepsilon_1^3))^2+\mathbb{E}(\varepsilon_1^6)$, we can bound the sum resulting from  this case of the
%%%% I MADE CHANGE HERE %%%%
$M=(\eta^2+\mathbb{E}(\varepsilon_1^6))$, we can bound the sum resulting from
%this case of the three by three configuration, as well as the one corresponding to all indices being equal by
these two configurations by
\begin{eqnarray*}\lefteqn{
\frac{M{q_2}^{-2d}}{n}\sum_{h=1}^{q_2}\sum_{h'=1}^{q_2}\sum_{t=0}^\infty\sum_{k=-n}^n\sum_{t'=0}^\infty a^2_ta_{t+k}a_{t'}a_{t'+h-k}a_{t'+h'}}\\
&&\le2\frac{M{q_2}^{-2d}}{n}\sum_{h=1}^{q_2}\sum_{h'=1}^{q_2}\sum_{k=0}^n\Delta(k)\gamma(h')\sim c\frac{{q_2}^{-2d}}{n}{q_2}n^d{q_2}^{2d}=o(n^{d-1/2})\to0.
\end{eqnarray*}
Consider now the configuration in (\ref{sumgamma}) where two indices are equal and another four are equal. An example of such configuration is when
$i-s=i-t\neq i+h-f=j-s'=j-t'=j+h'-f'$ so that $s=t$  and $s'=t'$, $f'=t'+h'$ and $f=t'+h+k$ where $k=i-j$. This configuration of the indices results in the following quantity
\begin{eqnarray*}
\frac{\mathbb{E}(\varepsilon_1^2)\mathbb{E}(\varepsilon_1^4)}{n}\sum_{h=1}^{q_2}\sum_{h'=1}^{q_2}\sum_{t=0}^\infty\sum_{t'=0}^\infty \sum_{k=-n}^n a_t^2a_{t'}^2a_{t'+h'}a_{t'+h+k},
\end{eqnarray*}
which is $O({q_2}^2/n)$.
This completes the proof of  (\ref{eq1}).\\

We now prove (\ref{eq2}). Noting that
$$
{q_3}^{-3d}\sum_{h=1}^{q_3}\sum_{h'=1}^{{q_3}-h}\left\vert\mathbb{E}\left(\hat\Delta(h,h')\right)-\Delta(h,h')\right\vert=O\left(\frac{{q_3}^2}{n^2}\right)\to0
$$
it will suffice to prove that
\begin{equation}\label{res}
{q_3}^{-6d}\mathbb{E}\left(\left(\sum_{h=1}^n\sum_{h'=1}^{{q_3}-h}\hat\Delta(h,h')\right)^2\right)-{q_3}^{-6d}\left(\sum_{h=1}^{q_3}\sum_{h=1}^{{q_3}-h}\Delta(h,h')\right)^2\to0.
\end{equation}
We have
\begin{eqnarray}\label{lastone}\lefteqn{
\mathbb{E}(\hat\Delta(h,h')\hat\Delta(u,u'))}\nonumber\\
&&=\frac{1}{n^2}\sum_{i=1}^{n-h-h'}\sum_{j=1}^{n-u-u'}
\sum_{s=0}^\infty
\sum_{t=0}^\infty\sum_{f=0}^\infty\sum_{s'=0}^\infty\sum_{t'=0}^\infty\sum_{f'=0}^\infty
a_aa_ta_fa_{s'}a_{t'}a_{f'}\mathbb{E}\left[\varepsilon_{i-t}\varepsilon_{i+h-s}\varepsilon_{i+h+h'-f}\varepsilon_{j-t'}
\varepsilon_{j+u-s'}\varepsilon_{j+u+u'-f'}\right].\nonumber\\
&&
\end{eqnarray}
We show that  the three by three configuration corresponding to $i-t=i+h-s=i+h+h'-f\neq j-t'=j+u-s'=j+u+u'-f'$ will coincide in the limit with
$$
{q_3}^{-6d}\left(\sum_{h=1}^{q_3}\sum_{h=1}^{{q_3}-h}\Delta(h,h')\right)^2.
$$
First, summing over $h,h',u,u',i,j$, this configuration of (\ref{lastone}) results in the following quantity (subtracting the diagonal elements with $k=j-i$)
\begin{eqnarray}\label{double1}\lefteqn{
\eta^2{q_3}^{-6d}\sum_{h=1}^{q_3}\sum_{h'=1}^{{q_3}-h}\sum_{u=1}^{q_3}\sum_{u'=1}^{{q_3}-u}\sum_{t=0}^\infty a_ta_{t+h}a_{t+h+h'}\sum_{t'=0}^\infty a_{t'}a_{t'+u}a_{t'+u+u'}}\nonumber\\
&&-\frac{\eta^2{q_3}^{-6d}}{n}\sum_{h=1}^{q_3}\sum_{h'=1}^{{q_3}-h}\sum_{u=1}^{q_3}\sum_{u'=1}^{{q_3}-u}\sum_{t=0}^\infty\sum_{k=-n}^n a_ta_{t+h}a_{t+h+h'}a_{t+k}a_{t+k+u}a_{t+k+u+u'}.
\end{eqnarray}
The second term in (\ref{double1}) is nonnegative and bounded by
$$
\frac{{q_3}^{-6d}}{n}\sum_{h=1}^{q_3}\sum_{h'=1}^{{q_3}-h}\sum_{u=1}^{q_3}\sum_{u'=1}^{{q_3}-u}\Delta(h,h')\Delta(u,u')=O(1/n)\to0.
$$
The first term in (\ref{double1})  equals
$$
\left({q_3}^{-3d}\sum_{h=1}^{q_3}\sum_{h'=1}^{{q_3}-h}\Delta(h,h')\right)^2.
$$
Therefore, in order to complete the proof of (\ref{res}), it suffices to show that all the other configurations of (\ref{lastone}) result in quantities converging to zero. Here again we will show this for selected cases of each remaining configuration. In the three by three configuration, when we take
$i-t=i+h-s=j-t'$ and $i+h+h'-f=j+u-s'=j+u+u'-f'$,  (note that this also includes the all-indices-are-equal configuration), with $k=j-i$, we see that, with $M=\mathbb{E}(\varepsilon_1^6)+\eta^2$. Such a configuration results in a quantity bounded by
\begin{eqnarray*}\lefteqn{
\sum_{h=1}^{q_3}\sum_{h'=1}^{q_3}\sum_{u=1}^{q_3}\sum_{u'=1}^{q_3}
\frac{M{q_3}^{-6d}}{n}\sum_{k=-n}^n\sum_{t=0}^\infty\sum_{f=0}^\infty a_ta_{t+h}a_{t+k}a_fa_{f+k+u-(h+h')}a_{f+k+u+u'-(h+h')}}\\
&&=\frac{M{q_3}^{-6d}}{n}\sum_{u=1}^{q_3}\sum_{u'=1}^{q_3}\sum_{h=1}^{q_3}\sum_{h'=1}^{q_3}\sum_{t=0}^\infty\sum_{k=0}^na_ta_{t+h}a_{t+k}\sum_{f=0}^\infty a_fa_{f+k+u-(h+h')}a_{f+k+u+u'-(h+h')}\\
&&+\frac{M{q_3}^{-6d}}{n}\sum_{u=1}^{q_3}\sum_{u'=1}^{q_3}\sum_{h=1}^{q_3}\sum_{h'=1}^{q_3}\sum_{t=0}^\infty\sum_{k=1}^na_ta_{t+h+k}a_{t+k}\sum_{f=0}^\infty a_fa_{f+k+u-(h+h')}a_{f+k+u+u'-(h+h')}\\
&&\le c\frac{M{q_3}^{-6d}}{n}\sum_{u=1}^{q_3}\sum_{u'=1}^{q_3}\sum_{h=1}^{q_3}\sum_{h'=1}^{q_3}n^d\gamma(h')\gamma(u')=
O\left(\frac{{q_3}^{-6d}}{n}{q_3}^{2+4d}n^d\right)=O\left(\frac{{q_3}^2}{n}\right)^{1-d}\to0.
\end{eqnarray*}
As to showing the negligibility of the  two by four configuration of (\ref{lastone}) we consider, for example,  the case when
the first two indices equal and the remaining four are equal. For this,  we get (with $k=i-j$) a term bounded by
\begin{eqnarray*}\lefteqn{
\frac{\mathbb{E}(\varepsilon_1^2)\mathbb{E}(\varepsilon_1^4){q_3}^{-6d}}{n}\sum_{h=1}^{q_3}\sum_{h'=1}^{q_3}\sum_{u=1}^{q_3}\sum_{u'=1}^{q_3}\sum_{t=0}^\infty\sum_{t'=0}^\infty\sum_{k=-n}^n
a_ta_{t+h}a_{t'+h+h'+k}a_{t'}a_{t'+u}a_{t'+u+u'}}\\
&&\le 2\frac{\mathbb{E}(\varepsilon_1^2)\mathbb{E}(\varepsilon_1^4){q_3}^{-6d}}{n}\sum_{h=1}^{q_3}\sum_{h'=1}^{q_3}\sum_{u=1}^{q_3}\sum_{u'=1}^{q_3}\sum_{t=0}^\infty\sum_{t'=0}^\infty\sum_{k=0}^n
a_ta_{t+h}a_ka_{t'}a_{t'+u}a_{t'+u+u'}\\
&&=2\frac{\mathbb{E}(\varepsilon_1^2)\mathbb{E}(\varepsilon_1^4){q_3}^{-6d}}{n}\sum_{h=1}^{q_3}\sum_{h'=1}^{q_3}\sum_{u=1}^{q_3}\sum_{u'=1}^{q_3}\sum_{k=0}^na_k\gamma(h)\Delta(u,u')
=O\left(\frac{{q_3}n^d{q_3}^{2d}{q_3}^{3d}}{{q_3}^{6d}n}\right)=O\left(\frac{{q_3}}{n}\right)^{1-d}\to0.
\end{eqnarray*}
This completes the proof of (\ref{res}) and that of (\ref{eq2}).
\bibliography{plain}

\begin{thebibliography}{99}
\bibitem {b} Abadir, M. K., Distaso, W., Giraitis, L. (2009). Two estimators for the long-run variance: Beyound short memory.Journal of Econometrics, 150. 56-70.
\bibitem{c} Bai, J. and  Ng S. (2005). Tests for Skewness, Kurtosis, and Normality for Time Series Data. Journal of Business and Economic Statistics. Vol. 23. 49-60.
\bibitem {z} Davydov, Y. (1970). The invariance principle for stationary processes (translated). Theor. Probability
App. 15, 487-498.
\bibitem{d} Grigoletto, M. and   Lisi, F. (2009). Looking for skewness in financial time series. Econometrics Journal. , volume 12, pp. 310-323.
\bibitem{Ko} Kokozska, P.s., and Taqqu, M.s. Stochastic Processes and their Applications 60 (1995) 19-47.
\end{thebibliography}

\end{document}